\documentclass[12pt,a4paper]{article}
\usepackage{authblk,amsmath,amsfonts,amssymb,amsthm,graphicx,caption,subcaption}
\usepackage[T2A]{fontenc}
\usepackage[utf8]{inputenc}
\usepackage[english]{babel}
\usepackage[top=2.0cm,bottom=2.0cm,left=2.0cm,right=2.0cm]{geometry}
\usepackage{multirow}
\providecommand{\keywords}[1]{\textbf{\textit{Keywords:}} #1}

\newcommand{\MM}{\mathbf{M}}
\newcommand{\KM}{\mathbf{K}}
\newcommand{\fM}{\mathbf{f}}
\newcommand{\aM}{\mathbf{a}}
\newcommand{\NM}{\mathbf{N}}
\newcommand{\BM}{\mathbf{B}}

\begin{document}

\begin{center}
\Huge{Numerical modeling of thin anisotropic composite membrane under dynamic load}\\[1cm]
\Large{V.V. Aksenov, A.V. Vasyukov, I.B. Petrov}\\[0cm]
\large{MIPT}\\[0cm]
\large{141700, Moscow region, Dolgoprudny, Institutsky lane, 9}\\[0cm]
\large{E-mail: aksenov.vv@phystech.edu}\\[0cm]
%\large{Поступила в редакцию: }\\[1cm]
\end{center}

\abstract{%
    This work aims to describe a mathematical model and a numerical method to simulate a thin anisotropic composite membrane moving and deforming in 3D space under a dynamic load of an arbitrary time and space profile. 
    The model and the method allow to consider problems when quasi static approximation is not valid, and elastic waves caused by the impact should be calculated. 
    The model and the method can be used for numerical study of different processes in thin composite layers, such as shock load, ultrasound propagation, non-destructive testing procedures, vibrations. 
    The thin membrane is considered as 2D object in 3D space, this approach allows to reduce computational time still having an arbitrary material rheology and load profile.
}

\keywords{composite membrane, dynamic load, numerical modeling}

%\fundnote{The work was supported by RFBR project 18-29-17027}

\bigskip

\section{Introduction}

Fabric composites are becoming an essential part of impact shields that protect spacecrafts from micrometeoroids and orbital debris. Present protection systems are based on the original Whipple shield \cite{whipple} that consists of two rigid (typically metal) layers spaced some distance from each other. This design allows to deal with particles with velocities up to 10-18 kilometres per second. An additional protection is provided by stuffed Whipple shields \cite{christiansen1995}, that include high-strength materials between rigid layers. The International Space Station uses different types of Whipple shields widely \cite {nasa}. Moreover, the International Space Station starts to use experimental inflatable modules, that use flexible fabric protection system without rigid layers that presented in the original Whipple shield. Computer design and optimization of stuffed Whipple shields and flexible protection system requires modeling of high-strength fabric composites under dynamic loading.

There is a number of works on numerical and experimental studies of fabric materials under a shock load, naming \cite{kobylkin_selivanov, walker1999, walker2001, porval} as an examples. Most of them consider a single relatively large impactor moving with the velocity 300-700 meters per second. However, for stuffed Whipple shields the load is completely different -- the initial particle is destroyed by the outer rigid layer, and the fabric layer is exposed to a cloud of small fragments distributed over time and space and moving with the velocity around 5-9 kilometres per second.

From mathematical model point of view this load profile means that it is not possible to consider the problem to be quasi static, since the speed of the fragments is comparable with the sound speed in the composite membrane. The stress-strain state of the composite layer should be calculated as the dynamic problem, and the elastic waves caused by the impact should be analyzed. Modeling of composite layer under asymmetric load requires also using anisotropic tensor of material elastic parameters and having three degrees of freedom for each point.

This work aims to provide all these features still describing composite layer as 2D object, since this allows to reduce computational time compared with 3D models. The same mathematical model and numerical method can be used later for numerical study of different processes in thin composite layers, such as ultrasound propagation, non-destructive testing procedures, vibrations.

This work is based on the models for thin fibers and membranes from \cite{rakhmatulin}. This work uses the same approach as \cite{rakhmatulin} to describe thin flexible structures, after that we do not obtain analytical solutions for particular cases, but solve the equations numerically and study a convergence rate of the numerical scheme for an arbitrary membrane and load.

\section{Mathematical model} 
\subsection{Notation}
    All the computations below are performed in Cartesian coordinates with $OX$, $OY$ axes lying in the plane of the undeformed membrane and $OZ$ axis is orthogonal to that plane. Vector $(x_0, y_0, z_0)^T$ denotes the initial coordinates of the point of the specimen in the undeformed state, $(x, y, z)^T$ coordinates of the same point at the given moment. Let 
$$
\mathbf{u} = \begin{bmatrix}
            u \\
            v \\
            w
            \end{bmatrix}
        = \begin{bmatrix}
            x - x_0 \\
            y - y_0 \\
            z - z_0
            \end{bmatrix}
$$
be the displacement vector
\subsection{Assumptions}
The membrane is considered to be an object with a characteristic size in one dimension ($z$) several orders of magnitude less than in other two. Due to this fact, only the displacement of the midsurface of the membrane in considered \cite{liu2013}. Based on this, we suggest that difference in displacement in $z$-direction can be neglected: 
\begin{equation}\label{no_z}
    \mathbf{u}(x, y, z) = \mathbf{u}(x, y)
\end{equation}
Thus, the membrane effectively is a 2D object in 3D space.

Only small strains are considered. Cauchy's strain tensor in the following form is used:
\begin{equation}
    \varepsilon_{ij} = \frac{1}{2}\left( \frac{\partial u_i}{\partial x_j} + \frac{\partial u_j}{\partial x_i} \right)
\end{equation}
We suppose that the materials are subject to Hooke's law in deformation rates considered:
\begin{equation}
    \sigma_{ij} = C_{ijkl}\varepsilon_{kl}
\end{equation}
The assumptions of small strains and linear elasticity up to destruction are generally valid for rigid composites during high speed interactions \cite{beklemysheva2016}. Other materials may show different behaviour, in this case they will not be covered by the model presented in this work.

\subsection{Stress-strain and strain-displacement relations}
We use the following vector notation for displacements:
    \begin{equation}\label{eq:deformation}
        \varepsilon = \begin{bmatrix}
                        \frac{\partial u}{\partial x} &
                        \frac{\partial v}{\partial y} &
                        \frac{\partial w}{\partial z} &
                        \frac{\partial u}{\partial y} + \frac{\partial v}{\partial x} &
                        \frac{\partial v}{\partial z} + \frac{\partial w}{\partial y} &
                        \frac{\partial w}{\partial x} + \frac{\partial u}{\partial z} 
                    \end{bmatrix}^T = \mathbf{Su}
    \end{equation}
    with linear differential operator $\mathbf{S}$ defined as
    \begin{equation}\label{eq:strain}
        \mathbf{S} =\begin{bmatrix}
                            \frac{\partial}{\partial x} &   0               &   0   \\  
                            0                           &   \frac{\partial}{\partial y} &   0   \\  
                            0       &   0               &   \frac{\partial}{\partial z} \\

                            \frac{\partial}{\partial y} &   \frac{\partial}{\partial x} &   0   \\  
                            0   &   \frac{\partial}{\partial z} &   \frac{\partial}{\partial y} \\  
                            \frac{\partial}{\partial z} &   0   &   \frac{\partial}{\partial x}
                    \end{bmatrix} 
    \end{equation}
The Hooke's law then takes the form:
    \begin{equation}\label{eq:stress}
        \mathbf{\sigma} = \begin{bmatrix} 
        \sigma_{xx} & 
        \sigma_{yy} & 
        \sigma_{zz} &
        \tau_{xy}   &
        \tau_{yz}   &
        \tau_{xz}
        \end{bmatrix}^T = \mathbf{D\varepsilon}
    \end{equation}
with $D$ being the compliance matrix. In the most general case $\mathbf{D}$ is symmetric and depends on $21$ independent elastic moduli. 
\subsection{Equations of motion}
We derive the equations of motion using virtual work principle. 
Let us first consider the equilibrium conditions for the unit volume under external loading in the static case. 
Let $\delta u, \delta \varepsilon$ be variations of corresponding values. $\mathbf{\overline{b}}$ denotes the distributed external force per unit volume.
Then the virtual work per unit volume is given by:
    \begin{equation}\label{eq:virtual_work}
        \delta w = \delta\mathbf{\varepsilon^T \sigma} - \delta\mathbf{u^T \bar{b}}
    \end{equation}
With $u$ sufficiently smooth $\delta\varepsilon = S\delta\mathbf{u}$.
For transition to dynamic case, according to d'Alembert's principle, we shall add distributed inertia force:
$$
\bar{\mathbf{b}} = \mathbf{b}  -\rho\ddot{\mathbf{u}},
$$
    where $\mathbf{b}$ is distributed load,  $-\rho\ddot{\mathbf{u}}$~---~inertia force.
Taking this into account along with \eqref{eq:strain}, \eqref{eq:stress}, we get the equation
$$
\delta\mathbf{u^T\left(S^TDSu + \rho\ddot{u} - b\right)} = 0,
$$
that must be satisfied at any variation $\delta\mathbf{u}$. Thus the term in brackets must be equal to zero. We finally get the equation of motion:
\begin{equation}\label{eq:motion}
    \mathbf{S^TDSu + \rho\ddot{u} - b} = 0
\end{equation}
The further expansion of this equation appears encumbering, in particular when $\mathbf{D}$ is given by more than two independent parameters, as for isotropic material. On the other hand, precisely those cases are of interest for modelling composite materials. Thus, in the following section, we shall derive the governing equations for the computational element. 
\section{Numerical method}
    Finite Element Method (FEM) was chosen for the task. 
One of the advantages of this approach is the possibility for treating complex geometries, laying down the basis for future study of the movement of perforated constructions or screens with holes created by the impacts. 
Formulation of the computational algorithm follows the general methodology described in \cite{zie00}. 
This paper is different from the traditional approach in the fact that material points are parametrized with two coordinates $(x_0, y_0)$, but has three degrees of freedom $(u, v, w)$. 
Traditionally, for 2D problems only formulations with either one (normal, $w$) or two (in-plane, $(u, v)$) degrees of freedom are studied. 
This is because for the isotropic material due to particular structure of matrix $\mathbf{D}$, the system of equations of motion is split into two systems:
wave equation for $w$ and Lam\'e equations for in-plane motion, which can be solved independently. 
This, though, is not true for composite materials, so the need for simultaneous solution of equation for all three degrees of freedom arises.
Compared with direct modelling with 3D FEM the proposed methods saves computational time, as we don't need to mesh throgh thickness direction $z$.

\subsection{Domain decomposition}
    
    The computational domain corresponds with the physical membrane. Unstructured grid of triangular elements is used. 
    We suppose thickness $h$, density $\rho$ and elastic moduli  $E, \nu$~of the material to be constant around the element. 
    The external force $\mathbf{b}$ is suggested to be distributed uniformly in the element. 
    The following values are associated with the verices: their initial coordinates $(x_0, y_0)$, displacements $\mathbf{a}$ and velocities $\mathbf{v}$.
    
    \begin{figure}[!h]
    \noindent\centering{
        \includegraphics[width=.7\textwidth]{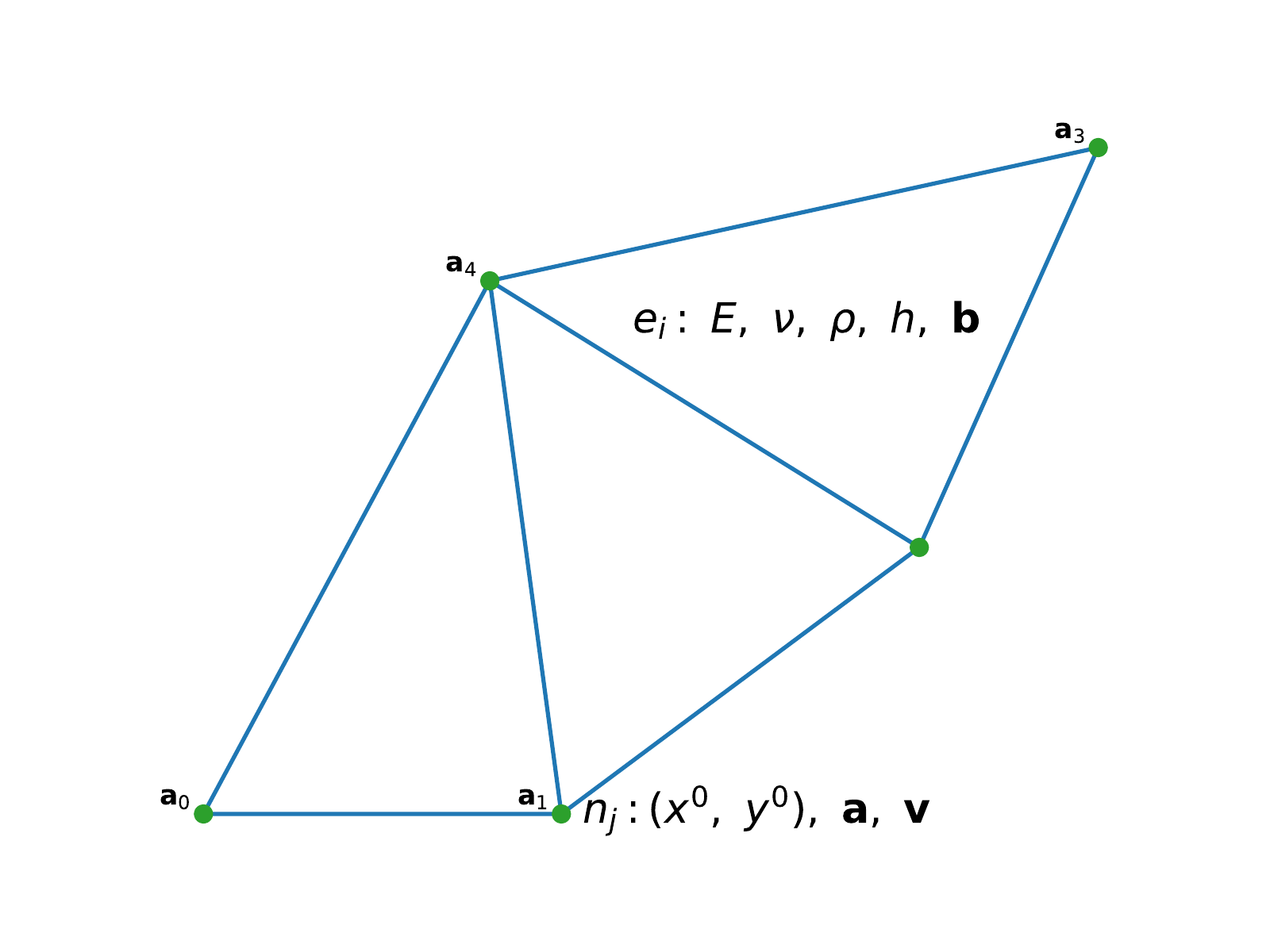}}
        \caption{Sample domain decompositon. Elements $e_i$ and vertices $n_j$ are presented along with corresponding values}
        \label{fig:dom_div}
    \end{figure}

\subsection{Displacement approximation}
    Consider a triangular element with vertices $(x_0^m, y_0^m), (x_0^n, y_0^n), (x_0^p, y_0^p)$. Let us denote the displacement vector at vertex $i$ as
    $$
        \mathbf{a}_i = \mathbf{u}(x_0^i, y_0^i)
    $$
    Displacements are approximated in the domain of the element as linear functions of nodal displacements by defining the \emph{shape functions} $\mathbf{N_i}(x, y, z)$:
        \begin{equation}\label{shape}
            \mathbf{u}(x,y,z) = \sum_{i \in \{m,n,p\}}{\mathbf{N_ia_i}} = \mathbf{Na}
        \end{equation}
    with $\NM = \left[ \begin{smallmatrix} \NM_m & \NM_n & \NM_p
                                            \end{smallmatrix} \right]$, 
        $\aM = \left[ \begin{smallmatrix} \aM_m \\ \aM_n \\ \aM_p
                        \end{smallmatrix} \right]$~---~displacements of the vertices,
    as only linear approximation is considered. Taking (\ref{no_z}) into account one gets:
        \begin{equation}\label{eq:shape_function}
            \NM_i = (\alpha_i + \beta_i x + \gamma_i y)\mathbf{E}
        \end{equation}
        with coefficients defined by:
        \begin{equation}
           \NM_i(x_j, y_j) = \begin{cases} 
                                \mathbf{E},\ i = j \\
                                \mathbf{0},\ i \neq j
                    \end{cases}\ i,\ j \in \{m,n,p\}
        \end{equation}
        Final coefficients for given shape functions are given by:
        \begin{gather}
            \alpha_i = \frac1{S_e}\begin{vmatrix}   
                                    x_j & y_j \\
                                    x_k & y_k  
                                    \end{vmatrix}    \quad
            \beta_i = -\frac1{S_e} (y_k - y_j) \quad
            \gamma_i = \frac1{S_e} (x_k - x_j) \quad
            S_e = \begin{vmatrix}
                        1   & x_i   & y_i \\
                        1   & x_j   & y_j \\
                        1   & x_k   & y_k
                    \end{vmatrix}\notag
        \end{gather}

\subsection{Approximate stress and strain }
    Let us approximate \eqref{eq:deformation} by substituting displacement with approximation from \eqref{shape} 
    \begin{equation} \label{eps}
                \mathbf{\varepsilon} = \mathbf{Su} = \mathbf{SNa} = \mathbf{Ba}
     \end{equation}
     Taking into account the formulas for $\NM$ the matrix $\BM$ becomes:
     \begin{equation}
        \mathbf{B} = \begin{bmatrix}
                \mathbf{B}_1^1 & \mathbf{B}_2^1 & \mathbf{B}_3^1 \\
                \mathbf{B}_1^2 & \mathbf{B}_2^2 & \mathbf{B}_3^2 
                \end{bmatrix}
    \end{equation}
    \begin{equation}
        \mathbf{B}_i^1 = \begin{bmatrix}
                            \beta_i & 0         & 0 \\
                            0       & \gamma_i  & 0 \\
                            0       & 0         & 0
                        \end{bmatrix}, \quad
        \mathbf{B}_i^2 = \begin{bmatrix}
                            \gamma_i    & \beta_i   & 0 \\
                            0       & 0             & \gamma_i\\
                            0       & 0             & \beta_i 
                        \end{bmatrix}
    \end{equation}
     Stress and strain are connected via Hooke's law \eqref{eq:stress}. Substituting with the equatiom for deformation from above, one gets:
    \begin{equation}\label{sigma}
       \sigma = \mathbf{D\varepsilon} = \mathbf{DBa}
    \end{equation}

    In the current implementation it is possible to define either all $21$ independent elastic constants, or, for the case of isotropic material, Young's modulus $E$ and Poisson's ratio $\nu$.
    The compliance matrix in the former case has the following form:
    \begin{gather}
       \mathbf{D} = 
        \begin{bmatrix}
            \mathbf{D}_1    & \mathbf{0}    \\
            \mathbf{0}  & \mathbf{D}_2
        \end{bmatrix}.
        \mathbf{D}_1 = \frac{E}{(1 + \nu)(1 - 2\nu)}
        \begin{bmatrix}
            1 - \nu     & \nu       & \nu \\
            \nu         & 1 - \nu   & \nu \\
            \nu     & \nu       & 1 - \nu
        \end{bmatrix} \\ \notag
        \mathbf{D}_2 = \frac{E}{(1 + \nu)(1 - 2\nu)}
        \begin{bmatrix}
            \frac{1-2\nu}2  & 0         & 0 \\
            0       & \frac{1-2\nu}2& 0 \\
            0       & 0     & \frac{1-2\nu}2
        \end{bmatrix}
    \end{gather}

\subsection{Virtual work for the element}\label{local}
    Let us get the approximation for \eqref{eq:virtual_work} by substituting the variations 
    with finite element approximations (\ref{shape}),~(\ref{eps}):
     \begin{equation}
        \delta\mathbf{u} = \NM\delta\aM, \quad \delta \mathbf{\varepsilon} = \BM\delta\aM
    \end{equation}

    Taking these into account along with (\ref{sigma}), the virtual work equation \eqref{eq:virtual_work} becomes:
    \begin{equation}\label{eq:elem_virt_work}
        \delta w = \delta\mathbf{\varepsilon^T \sigma} - \delta\mathbf{u^T \bar{b}}
                 = \mathbf{\delta a^T\left(B^T D Ba  - N^T \bar{b}\right)}
    \end{equation}
    
    Let us integrate the equation over the element volume and introduce fictive nodal forces $\mathbf{q_i}$, that balance internal elastic forces and external loads:
    \begin{equation}\label{eq:elem_virt_work_integral}
        \delta \aM^T  \left(\int_{V_e}\mathbf{B^T D Ba}dV
                    - \int_{V_e} \mathbf{N^T\bar{b}}dV \right) = \delta\aM^T\mathbf{q}_e
    \end{equation}
    Transition to the dynamic case is done according to the D'Alembert principle:
    \begin{equation}\label{eq:inertion}
        \bar{\mathbf{b}} = \mathbf{b} - \rho \ddot{\mathbf{u}}
    \end{equation}
    where $\mathbf{b}$ is external force,  $-\rho\ddot{\mathbf{u}}$~---~force of inertia.
    For acceleration we shall use the same approximation as for displacement:
    \begin{equation}\label{eq:acceleration}
        \ddot{\mathbf{u}}(x, y) = \mathbf{N}(x, y)\ddot{\aM}. 
    \end{equation}
    Taking into account\eqref{eq:inertion} and \eqref{eq:acceleration} for computational element from \eqref{eq:elem_virt_work_integral} one finally gets:
    \begin{equation}\label{eq:single_elem_movement}
        \MM_e\ddot{\aM} + \KM_e\aM + \fM_e = \mathbf{q}_e
   \end{equation}
   Here the following matrices have been intriduced
   \begin{gather}
        \KM_e = \int_{V_e}{\mathbf{B^TDB}}dV
            \text{~---~stiffnes matrix}\\
        \MM_e = \rho\int_{V_e}{\mathbf{N^TN}}dV
            \text{~---~mass matrix}\\
        \fM_e = - \int_{V_e} \mathbf{N^Tb}dV_e 
            \text{~---~load vector}
   \end{gather}

\subsection{Global equations assembly}
    We now have to account for the effect from all elements that a certain vertex belongs to. 
    Instead of local displacement vector, defined in \ref{shape}, consider a $3N_{nodes}$-dimensional vector, consisting of stacked displacement vectors of each vertex.

    Thus we introduce global matrices $\KM, \MM$ with dimensions $3N_{nodes}\times3N_{nodes}$ and vector $\fM$ with dimension $3N_{nodes}$.
    They are assembled by following rules:
    \begin{gather}
        \KM_{ij} = \sum_{e}\KM^e_{ij} \\
        \MM_{ij} = \sum_{e}\MM^e_{ij} \\
        \fM_i = \sum_{e}\fM^e_{i}
    \end{gather}
    Here $\KM_{ij}$ denotes a $3\times 3$ block of $\KM$, standing at the intesections of row and column corresponding to $i$-th and $j$-th vertex. 
    $\KM_{ij}^e$ is a local stiffness matrix defined above in section \ref{local}. 
    Summation is done for all elements containing both i-th and j-th vertex.
    Fictive nodal forces are summed to, and, obviously, the sum equals to zero in case of equilibrium.
    Eventually one arrives at the following equation for the whole domain, analogous to \eqref{eq:single_elem_movement}:
    \begin{equation}\label{eq:movement}
        \boxed{\MM\ddot{\mathbf{a}} + \KM\aM + \fM = 0}
    \end{equation}

\subsection{Applying constraints}\label{sec:constraints}
    To study pinpoint strikes we need to constrain the velocity of the strike point. 
    A problem of studying the membrane with fixed border may also arise. 
    We shall demonstrate that such constraints can be taken into account without changing the structure of the governing equation (\ref{eq:movement}).
    Let us fix node $i$.
    In $\KM$ we fill the line $\KM_{ik}, k\in \overline{1, N_{nodes}}$ with zeros. 
    In $\MM$ let $M_{ii} = \mathbf{E}$, and the rest of the blocks $M_{ik}, k \neq i$ filled with zeros. Finally, in $\fM$ vector we fill $\fM_i$ with zeros.
    As a result, the system now contains an equation $\ddot{\mathbf{a}}_i = \mathbf{0}$, т.е. $\mathbf{v}_i(t) = \mathbf{v}_{fix}$. 
    To constrain displacement we shall determine $\mathbf{v}_{fix} = \mathbf{0}$
\subsection{Time integration}
    The procedure above allowed us to reduce an initial problem for PDE to an ODE problem \eqref{eq:movement}. 
    We then integrate this equation numerically using Newmark method \cite{newmark1959method}: 
    \begin{gather}
        \dot{\overline{\aM}}_n = \dot{\aM}_n + \tau(1-\beta_1)\ddot{\aM}_n \notag\\
        \overline{\aM}_n = \aM_n + \tau\dot{\aM}_n 
                + \frac12\tau^2(1 - \beta_2)\ddot{\aM}_n \notag\\
        \ddot{\aM}_{n+1} = - A^{-1}\left(\fM_{n+1} +\KM\overline{\aM}_n
                        \right),\quad A = \MM + \frac{1}{2}\tau^2\beta_2\KM\\
        \dot{\aM}_{n+1} = \dot{\overline{\aM}}_n + \beta_1\tau\ddot{\aM}_{n+1} \notag\\
        \aM_{n+1} = \overline{\aM} + \frac{1}{2}\tau^2\beta_2\ddot{\aM}_{n+1} \notag
    \end{gather}
    If $\beta_2 \geq \beta_1 \geq \frac12$, the method is unconditionally stable. 
    If $\beta_1 = \frac12$, then the method has the second order of approximation \cite{zie00}. 
    For all the numerical experiments below, parameters $\beta_1 = \beta_2 = \frac12$ were used.

\subsection{Numerical results}

\begin{figure}
    \noindent\centering{\includegraphics[width=.4\textwidth]{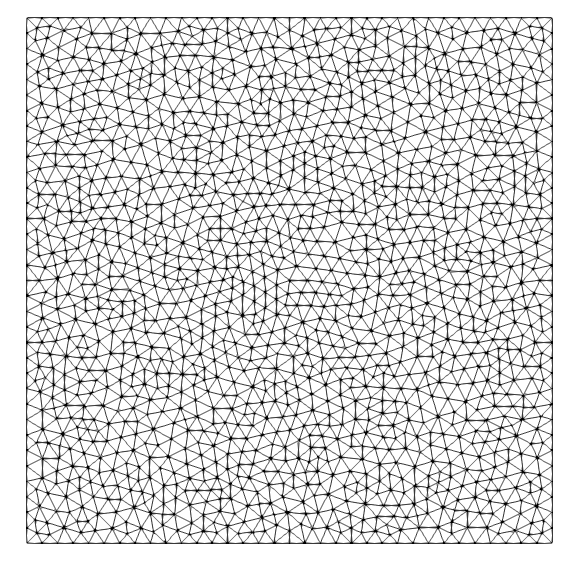}}
        \caption{Example calculation grid}
        \label{fig:non_uniform_mesh_sample}
    \end{figure}

The method described works on arbitrary irregular grids. The figure \ref{fig:non_uniform_mesh_sample} shows an example of a grid constructed using the gmsh mesh generator \cite{gmsh}.

\begin{figure}
  \begin{subfigure}{.45\linewidth}
    \centering\includegraphics[width=1.0\linewidth]{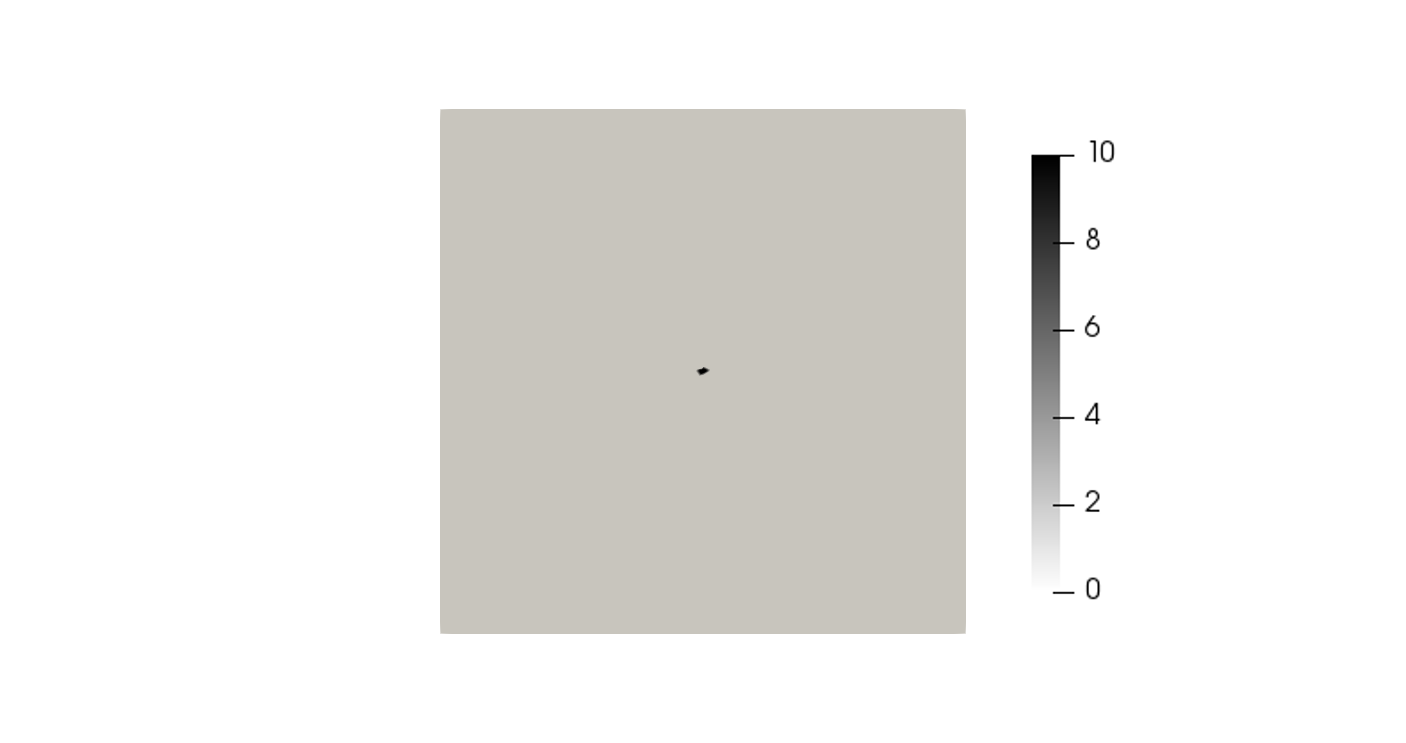}
    \caption{Initial state}
  \end{subfigure}
  \begin{subfigure}{.45\linewidth}
    \centering\includegraphics[width=1.0\linewidth]{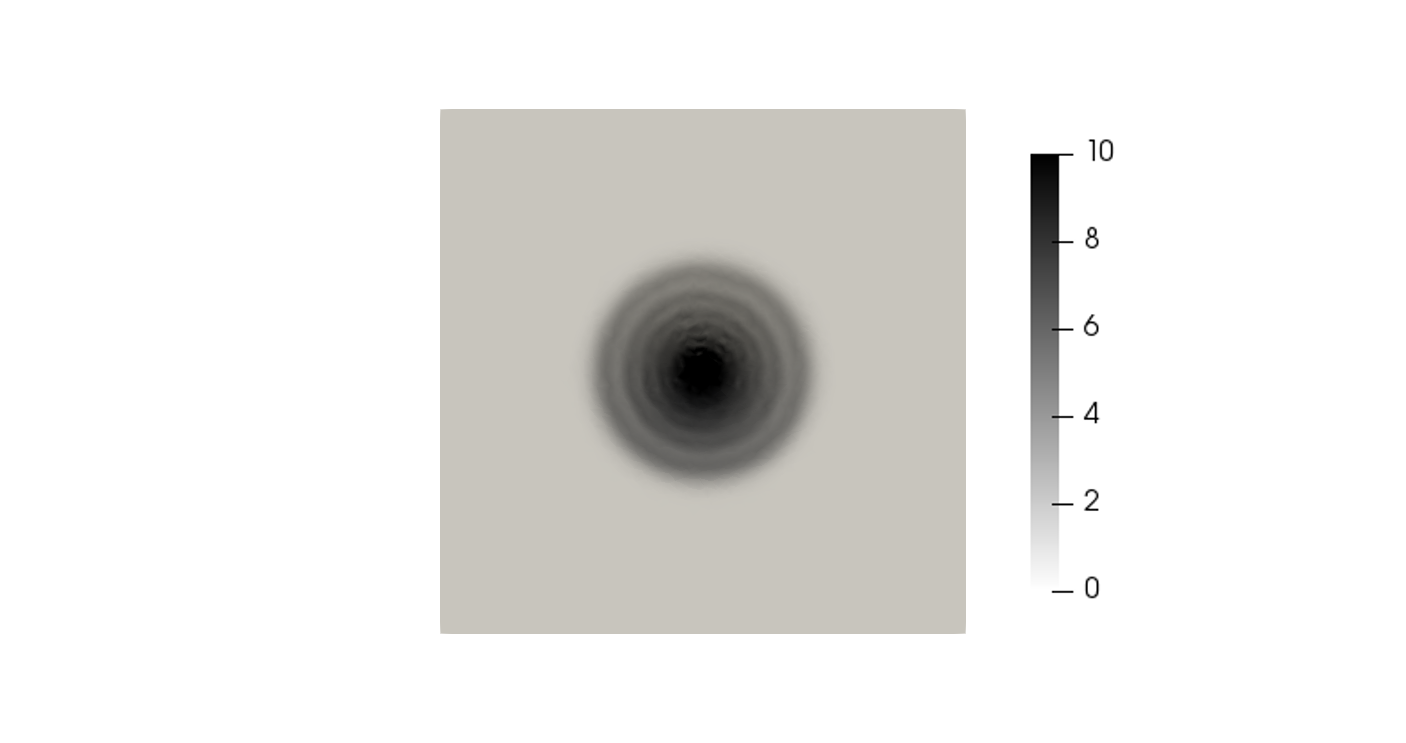}
    \caption{Wave propagation in the membrane}
  \end{subfigure}
 
  \begin{subfigure}{.45\linewidth}
    \centering\includegraphics[width=1.0\linewidth]{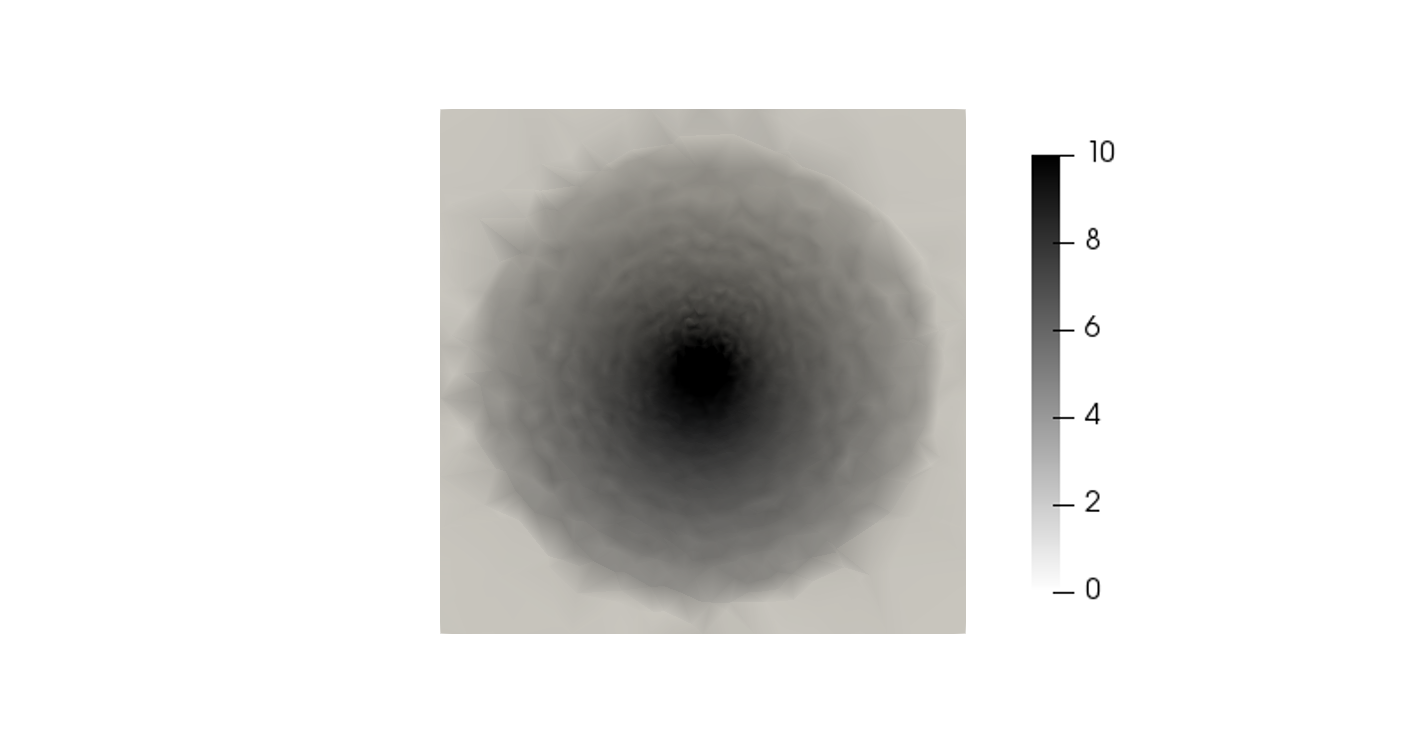}
    \caption{Wave reaches the border of the membrane}
  \end{subfigure}
  \begin{subfigure}{.45\linewidth}
    \centering\includegraphics[width=1.0\linewidth]{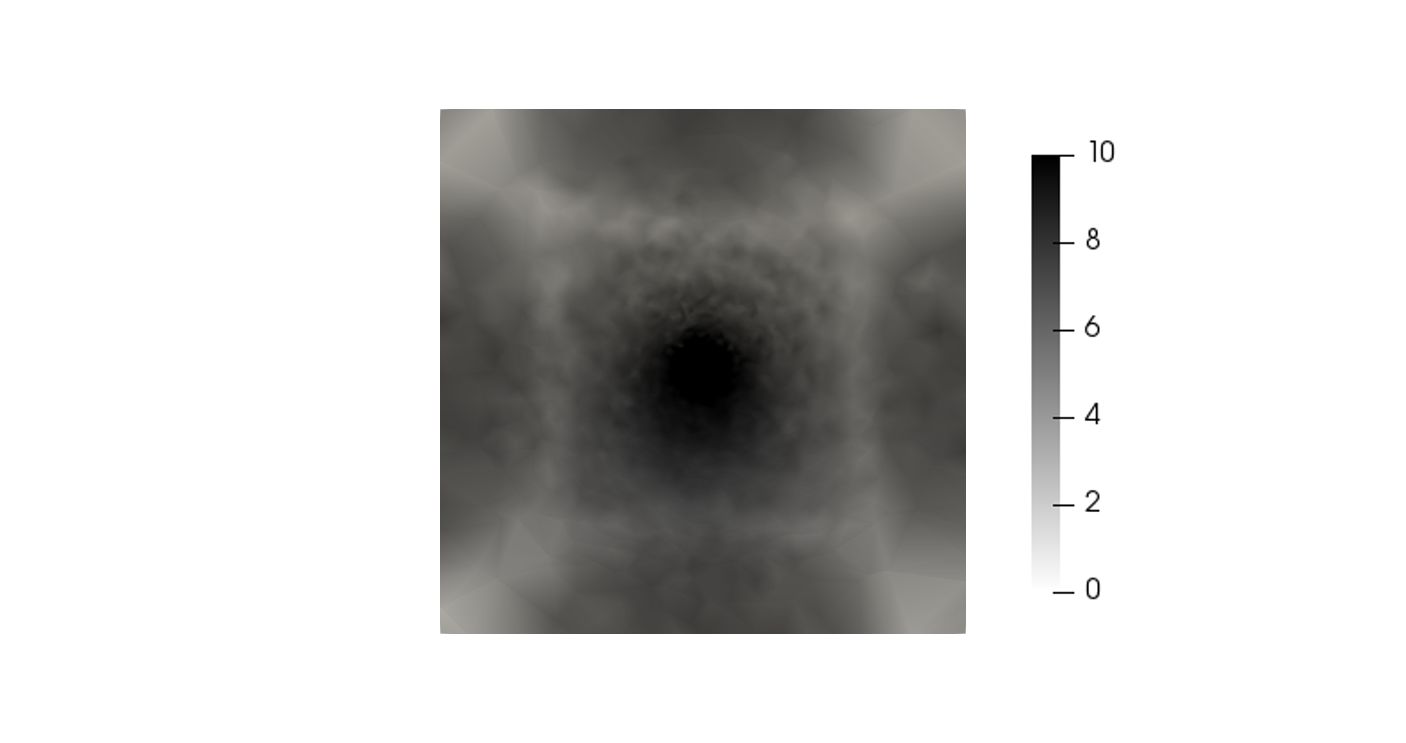}
    \caption{Wave is reflected from the border}
  \end{subfigure}

  \caption{Exposure to a point load impulse. The dynamics of the involvement of the membrane material in motion. The color shows the velocity module at different points in time.}
  \label{fig:normal_strike_waves}
\end{figure}

The figure \ref{fig:normal_strike_waves} shows an example of calculation using this grid. A single central mesh element is impacted at a constant speed directed normal to the membrane plane. The involvement of membrane material in movement is shown. It can be seen that for a symmetric formulation the solution is symmetric and does not contain numerical artifacts.

\begin{figure}
  \begin{subfigure}{.45\linewidth}
    \centering\includegraphics[width=1.0\linewidth]{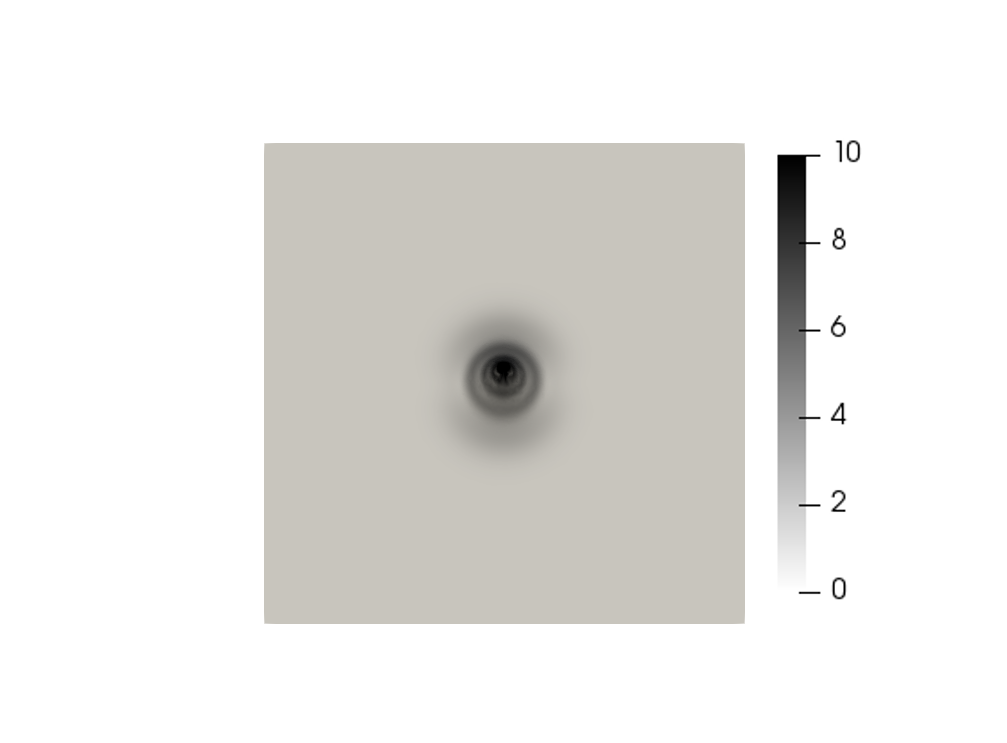}
    \caption{Initial stage of loading}
  \end{subfigure}
  \begin{subfigure}{.45\linewidth}
    \centering\includegraphics[width=1.0\linewidth]{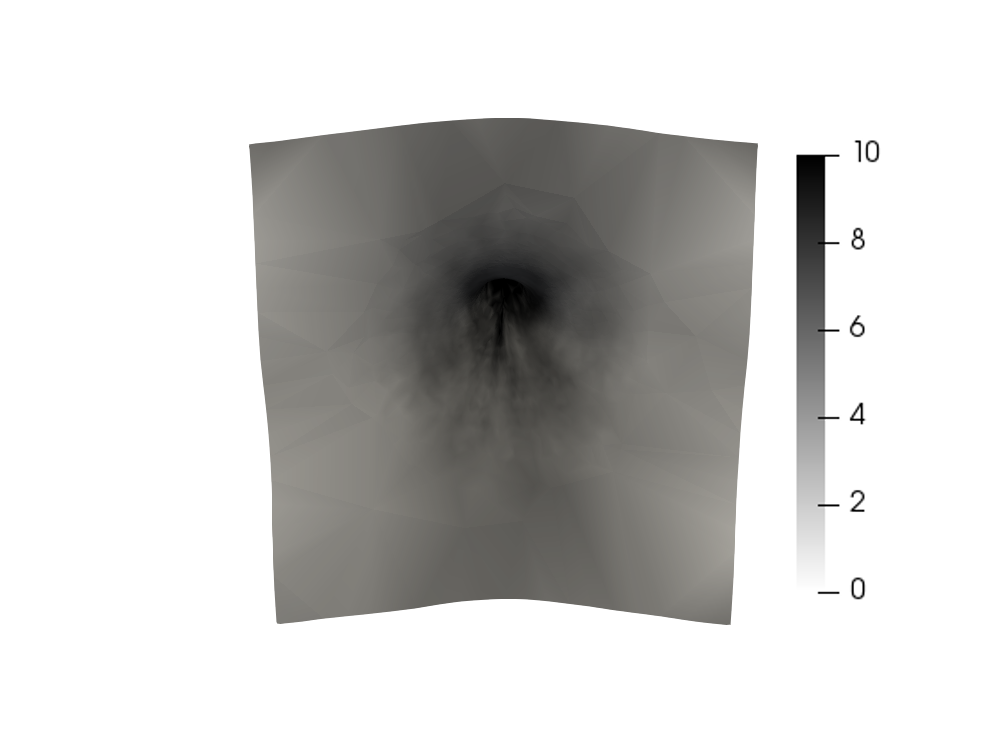}
    \caption{Asymmetric deformation}
  \end{subfigure}

  \caption{Exposure to a point load impulse at an angle of 30 degrees to the normal. The color shows the velocity module at different points in time.}
  \label{fig:skew_strike_waves}
\end{figure}

The implemented method allows to model asymmetric load profiles. The figure \ref{fig:skew_strike_waves} shows an example of a calculation in which a blow with a constant speed is applied at an angle to the normal. The involvement of the membrane material in motion and its substantially asymmetric deformation are seen.

\begin{figure}
  \begin{subfigure}{.45\linewidth}
    \centering\includegraphics[width=1.0\linewidth]{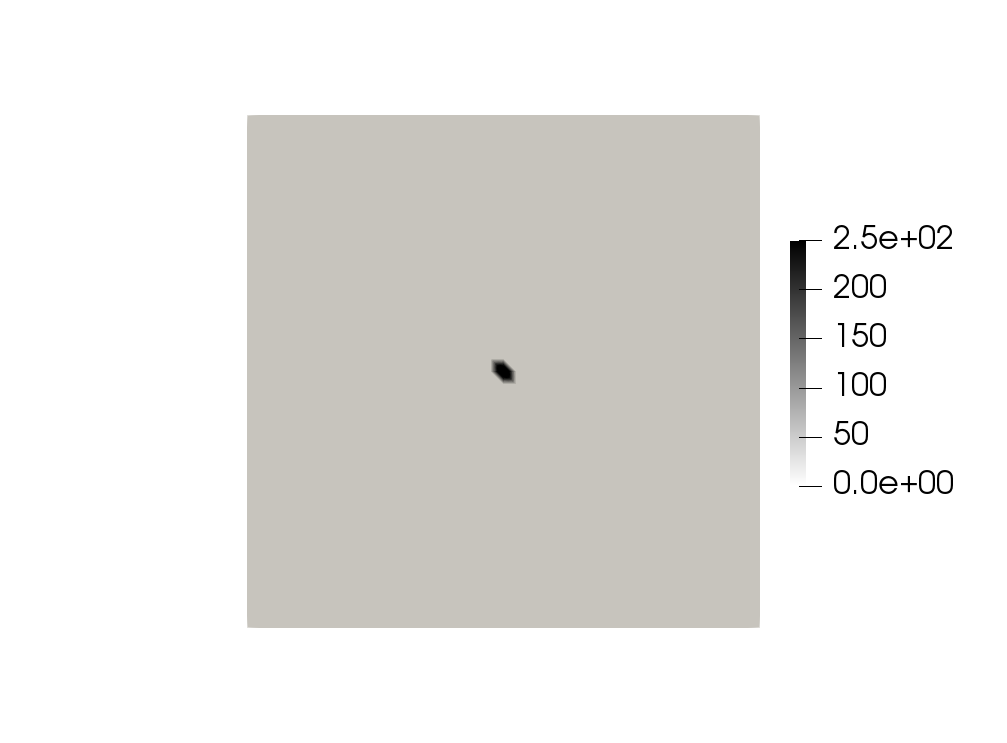}
    \caption{Initial stage of loading}
  \end{subfigure}
  \begin{subfigure}{.45\linewidth}
    \centering\includegraphics[width=1.0\linewidth]{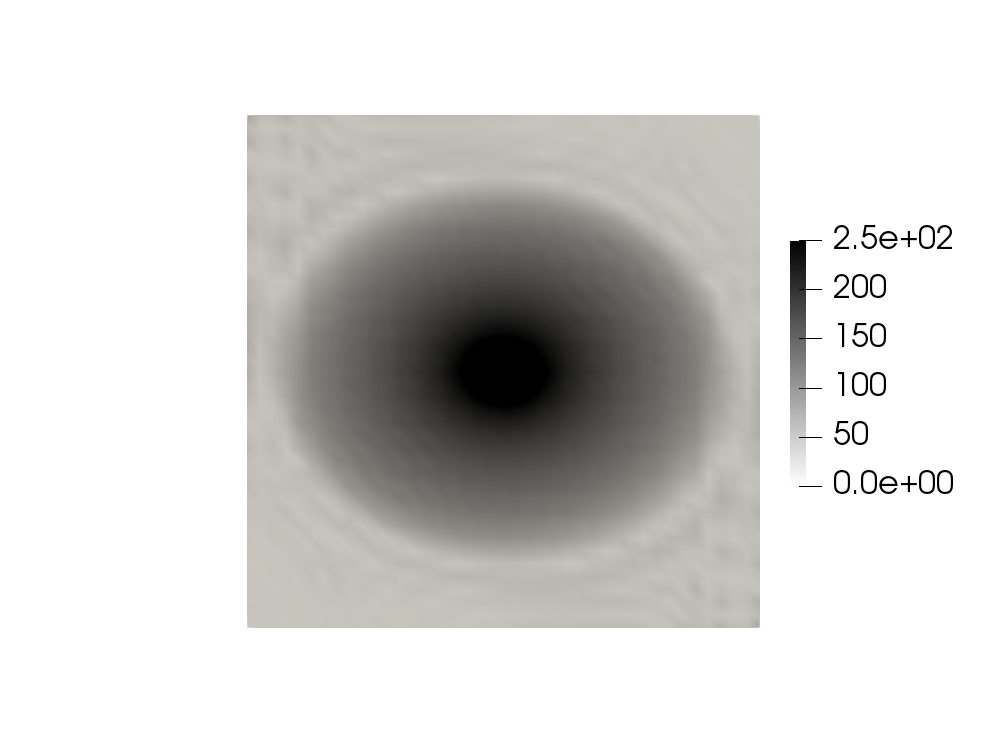}
    \caption{Wave reaches the border of the membrane}
  \end{subfigure}

  \caption{Anisotropic material exposured to a point load impulse. The color shows the velocity module at different points in time. Elliptic wave pattern formed due to material anisotropy.}
  \label{fig:aniso_material_waves}
\end{figure}

The figure \ref{fig:aniso_material_waves} shows an example of a calculation using anisotropic material model. A single central mesh element is impacted at a constant speed directed normal to the membrane plane. Elastic properies of the material are presented in table \ref{tab:material}. Elliptic wave pattern formed due to material anisotropy is seen of the figure.

\begin{table}[h]
    \centering
        \begin{tabular}[]{|c|l|}
            \hline
            $c_{11}$, GPa              &          150 \\
            $c_{12}$, GPa              &           40 \\
            $c_{13}$, GPa              &           10 \\
            $c_{22}$, GPa              &          150 \\
            $c_{23}$, GPa              &           80 \\
            $c_{33}$, GPa              &          150 \\
            $c_{44}$, GPa              &           80 \\
            $c_{55}$, GPa              &           20 \\
            $c_{66}$, GPa              &           30 \\
            \hline
        \end{tabular}
    \caption{Anisotropic material elasticity matrix non-zero elements}
    \label{tab:material}
\end{table}

\begin{figure}
    \noindent\centering{\includegraphics[width=.45\textwidth]{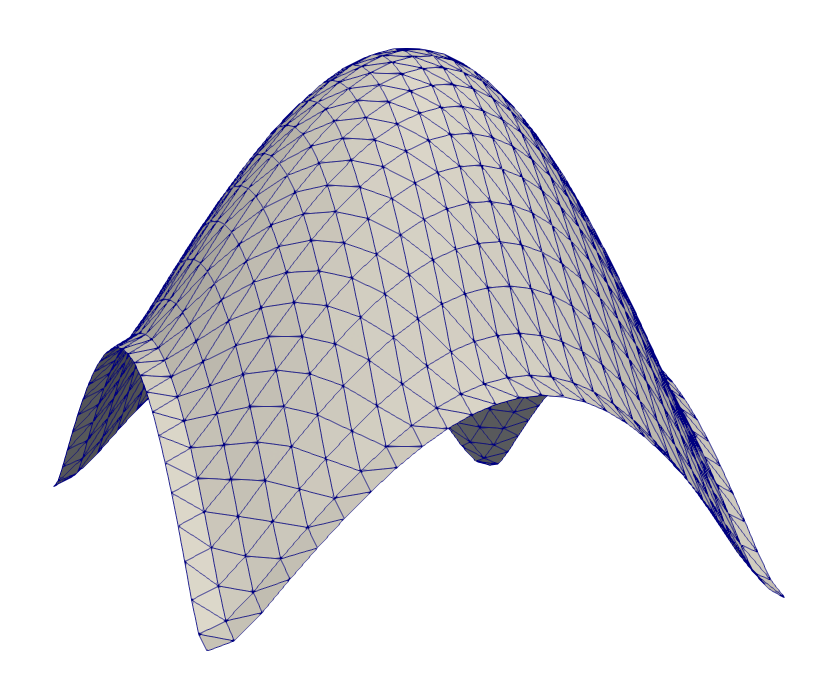}}
        \caption{Distributed load example}
        \label{fig:large_deformation_sample}
    \end{figure}

The method allows not only to obtain the initial membrane dynamics and elastic waves from the shock load, but also can be used to calculate of deformations that develop over a considerable time. For example, the figure \ref{fig:large_deformation_sample} shows the result of calculating the deformation of the membrane under the action of a distributed load, which is directed normal to the undeformed membrane and has the modulus
\begin{equation}\label{distrib_force}
    b(r) =  \begin{cases}
                b_0\cos^2{r}, \quad r \leq L \\
                0,              \quad \text{ otherwise}
            \end{cases} \\
            r = \frac{\pi}{2L}\sqrt{\left(x - \frac{1}{2}L \right)^2 
                            + \left(y - \frac{1}{2}L \right)^2},\
                            L\text{~---~membrane size}
\end{equation}

%\clearpage

\section{Actual convergence rate}
Scheme's actual convergence rate was examined with a method similar to presented in \cite{khokhlov2019}. 
First, a rather coarse regular grid, composed of rectangles, each split in two triangles by a diagonal, is generated. 
The choice of regular grid is dictated by the need to easily control the coarseness of the grid.
The solution is calculated up to some fixed time $T$. 
Velocities and displacements in all the points of this coarse grid are saved as a baseline solution. 
Then the calculations are repeated on a grid refined by splitting each rectangle into four smaller ones. Timestep is also scaled down accordingly.
The norms of difference of displacements and velocities for consequetive solutions on finer grids are calculated. 
Displacements and velocities for norm calculation are taken at the same space and time points where they were taken during the baseline solution.

Several test cases were examined, their descriptions summarized in Table \ref{tab:experiments}. 
In cases $1, 2$ the loading force is applied to one central square (a pair of triangular elements) in the middle of the membrane.
In cases $3, 4$ the velocity in the point closest to the geometric center of the membrane is constrained, the velocity vector being respectievely normal and inclined at an angle $\frac{\pi}6$ to normal.
\begin{table}[h]
    \centering
        \begin{tabular}[]{|c|l|}
            \hline
            Test case& Description\\
            \hline
            $1$ &  Normal strike with fixed load\\
            $2$ &  Strike with fixed load, inclined $\frac{\pi}{6}$  to normal\\
            $3$ &  Normal strike with fixed speed\\
            $4$ &  Strike with fixed speed, inclined $\frac{\pi}{6}$ to normal\\
            $5$ &  Normal load defined by formula \eqref{distrib_force} \\
            \hline
        \end{tabular}
    \caption{Test case formulations}
    \label{tab:experiments}
\end{table}

To demonstrate the convergence rate study, we present the plot used in convergence rate assessment for test case $3$.
The number of grid refinements $k$ is plotted over $OX$ axis, over $OY$ is plotted the norm of difference between solutions on consequtive grids in logarithmic scale.
Three norms: $L_{1}, L_{2}, L_{inf}$ are considered.
The slope of the fitted line shall represent the actual convergence rate. 

\begin{figure}
    \noindent\centering{        
    \includegraphics[width=.8\textwidth]{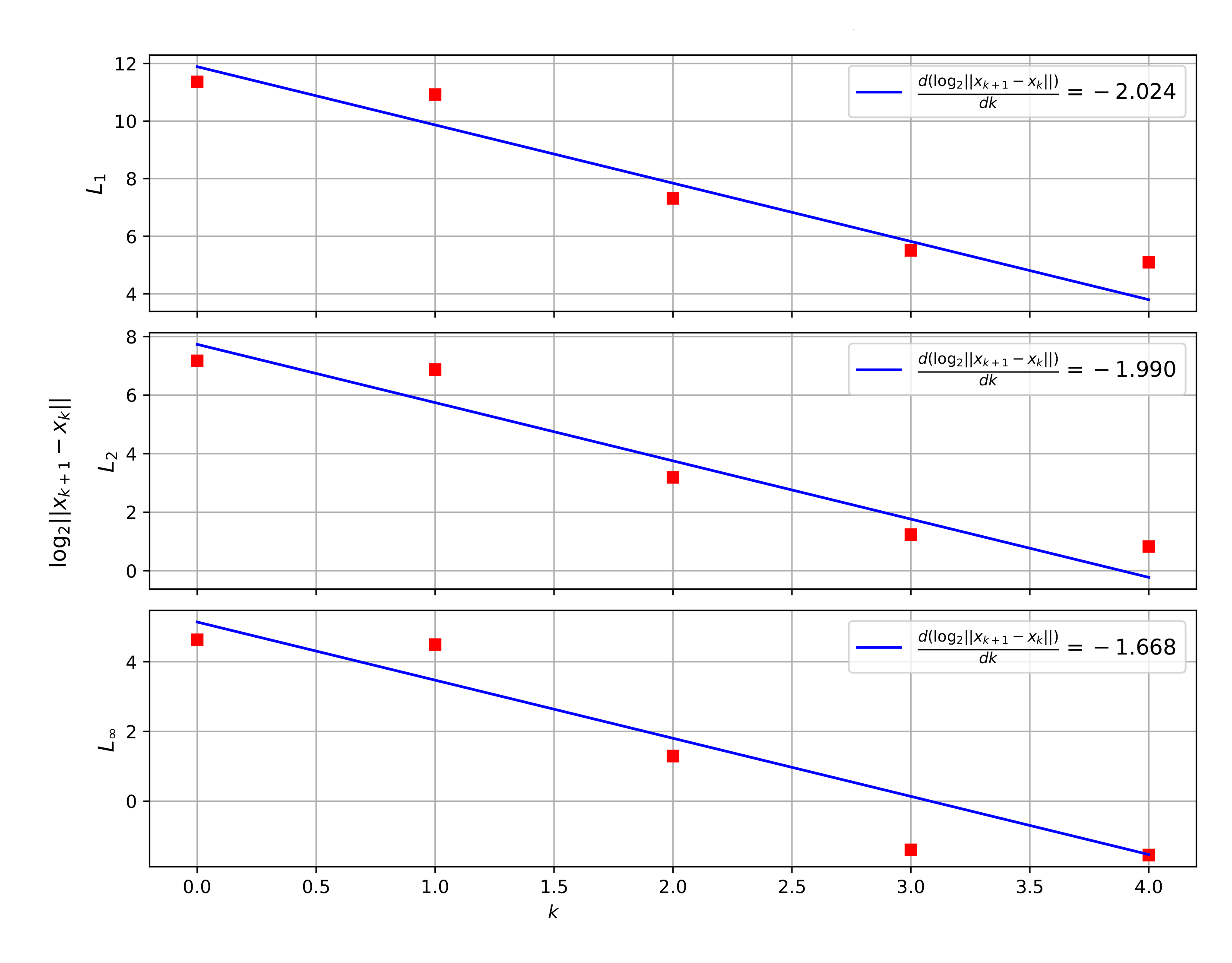}}
    \caption{Convergence: normal strike}
        \label{fig:conv_str}
\end{figure}

The results of convergence rate study are presented in \ref{tab:conv_results}

\begin{table}[h]
\centering
\begin{tabular}{|l|l|l|l|}
\hline
\multirow{2}{*}{Case} & \multicolumn{3}{c|}{Actual rate} \\ \cline{2-4} 
                            & $L_1$            & $L_2$            & $L_\infty$          \\ \hline
    $1$   & 2.558          & 2.504          & 2.511             \\ \hline
    $2$   & 1.534          & 1.482          & 1.445             \\ \hline
    $3$   & 2.024          & 1.990          & 1.668             \\ \hline
    $4$   & 1.599          & 1.660          & 1.964             \\ \hline
    $5$   & 1.748          & 1.678          & 1.561             \\ \hline
    \end{tabular}
\caption{Convergence rates}
\label{tab:conv_results}
\end{table}
%\clearpage
    
\section{Conclusions}

The paper describes the mathematical model and the numerical method for modeling a thin anisotropic composite membrane under a dynamic shock load of an arbitrary time and space profile. The results presented show that the convergence rate of the numerical scheme in the worst case is 1.4 (oblique strike when using the norm $L_{inf}$), and in the best case it is 2.5 (symmetric load using any norm). This order allows using the described numerical method for applied calculations carried out in conjunction with field experiments. Anisotropic materials supported by the numerical method can be used to describe multilayer fabric membranes using their effective parameters.

The model and the method are designed to consider the composite membrane as 2D object in 3D space still having an arbitrary material rheology and load profile, this approach allows to reduce computational time compared with direct modelling using 3D solvers.

An obvious limitation of the described model and method is the lack of the possibility of a detailed calculation of problems associated with a destruction of the membrane. Using the presented approach, it is only possible to determine a start of a failure using strain or stress thresholds. Enhacing the method to cover the destruction is the topic of future work.

\section*{Acknowledgements}

The work was supported by RFBR project 18-29-17027.

The authors are grateful to Beklemysheva K.A., Ph.D. for long thoughtful discussions of the results and Bot.Cafe for a warm atmosphere and great coffee.


\begin{thebibliography}{99}

\bibitem{whipple} Whipple F.L. Meteorites and Space Travel // Astronomical Journal, 1947, vol. 52, p.131.

\bibitem{christiansen1995} Christiansen E.L., Crews J.L., Williamsen J.E., Robinson J.H., Nolen A.M. Enhanced meteoroid and orbital debris shielding // International Journal of Impact Engineering, 1995, vol. 17, issues 1-3, pp. 217–228.

\bibitem{nasa} E. L. Christiansen, J. Arnold, A. Davis, J. Hyde et al. Handbook for Designing MMOD Protection // NASA Johnson Space Center, Houston, 2009.

%\bibitem{ahmad1970} Ahmad S., Irons B.M., Zienkiewicz O.C. Analysis of thick and thin shell structures by curved finite elements // International Journal for Numerical Methods in Engineering, 1970, vol. 2, issue 3, pp. 419-451.

%\bibitem{jenkins1996} Jenkins C.H. Nonlinear dynamic response of membranes: State of the art - Update // Applied Mechanics Reviews, 1996, vol. 49, issue 10S, pp. 41-48.

%\bibitem{sathe2007} Sathe S., Benney R., Charles R., Doucette E., Miletti J., Senga M., Stein K., Tezduyar T.E. Fluid–structure interaction modeling of complex parachute designs with the space–time finite element techniques // Computers \& Fluids, 2007, vol. 36, issue 1, pp. 127-135.

%\bibitem{miyazaki2006} Miyazaki Y. Wrinkle/slack model and finite element dynamics of membrane // International Journal for Numerical Methods in Engineering, 2006, vol. 66, issue 7, pp. 1179-1209.

\bibitem{kobylkin_selivanov} Kobylkin I.F., Selivanov V.V. Materials and Structures of Light Armor Protection // BMSTU, Moscow, Russia, 2014. (in Russian)

\bibitem{walker1999} Walker J.D. Constitutive Model for Fabrics with Explicit Static Solution and Ballistic Limit // Proceedings of the Eighteenth International Symposium on Ballistics, San Antonio, USA, 1999.

\bibitem{walker2001} Walker J.D. Ballistic Limit of Fabrics with Resin // Proceedings of the Nineteenth International Symposium on Ballistics, Interlaken, Switzerland, 2001.

\bibitem{porval} Porval P.K., Phoenix S.L. Modeling System Effects in Ballistic Impact into multi-Layered Fibrous Materials for Soft Body armor // International Journal of Fracture, 2005, vol. 135, issue 1–4, pp. 217-249.

%\bibitem{dmitrochenko2003} Dmitrochenko O.N., Pogorelov D.Y.U. Generalization of plate finite elements for absolute nodal coordinate formulation // Multibody System Dynamics, 2003, vol. 10, issue 1, pp. 17-43.

%\bibitem{shabana2003} Mikkola A.M., Shabana A.A. A non-incremental finite element procedure for the analysis of large deformation of plates and shells in mechanical system applications // Multibody System Dynamics, 2003, vol. 9, issue 3, pp. 283–309.

\bibitem{rakhmatulin} Rakhmatulin K.A., Demianov Y.A. Strength under high transient loads // Israel Program for Scientific Translations, 1966.

\bibitem{liu2013} Liu C., Tian Q,, Yan D. and Hu H. Dynamic analysis of membrane systems undergoing overall motions, large deformations and wrinkles via thin shell elements of ANCF // Computer Methods in Applied Mechanics and Engineering, 2013, vol. 258, pp. 81-95.

\bibitem{beklemysheva2016} Beklemysheva K.A., Vasyukov A.V., Ermakov A.S., Petrov I.B. Numerical simulation of the failure of composite materials by using the grid-characteristic method // Mathematical Models and Computer Simulations, 2016, vol. 8, issue 5, pp. 557-567.

\bibitem{zie00} Zienkiewicz O.C. and Taylor R.L. Finite Element Method: Volume 1 - The Basis // 5th Oxford: Butterworth-Heinemann, 2000.

\bibitem{newmark1959method} Newmark N.M. A Method of Computation for Structural Dynamics // Journal of the Engineering Mechanics Division, 1959, vol. 85, issue 3, pp. 67-94.

\bibitem{gmsh} Geuzaine C., Remacle J.-F. Gmsh: A 3-D finite element mesh generator with built-in pre- and post-processing facilities // International Journal for Numerical Methods in Engineering, 2009, vol. 79, issue 11, pp. 1309-1331.

\bibitem{khokhlov2019} Khokhlov N.I., Golubev V.I. On the Class of Compact Grid-Characteristic Schemes // In: Smart Modeling for Engineering Systems, 2019, pp. 64-77.

\end{thebibliography}
\end{document}